\input amstex

\documentstyle{amsppt}
\magnification=\magstep1 \pagewidth{5.0in} \pageheight{6.7in}
%\hcorrection{-0.4in}
%\vcorrection{-0.4in}
\abovedisplayskip=10pt \belowdisplayskip=10pt

\topmatter

\title New $q$-Euler numbers and polynomials associated with $p$-adic $q$-integrals
\endtitle

%thanks will be a 1st page footnote
\author{ Taekyun Kim, Min-Soo Kim, Leechae Jang and Seog-Hoon Rim  }\endauthor
\address{
\endgraf{EECS, Kyungpook National University, Taegu 702-701, South Korea}
\endgraf {\it Email Address}: {\rm tkim$\@$knu.ac.kr}
\endgraf{\;}
\endgraf{\;}
\endgraf{Department of Mathematics, Kyungnam University, Masan 631-701, South Korea}
\endgraf {\it Email Address}: {\rm mskim$\@$kyungnam.ac.kr}
\endgraf{\;}
\endgraf{\;}
\endgraf{Department \, of \, Mathematics \, and \, Computer Science,\, KonKuk University,\,
Chungju \endgraf 380-701, South Korea}
\endgraf {\it Email Address}: {\rm leechae.jang$\@$kku.ac.kr, \rm leechae-jang$\@$hanmail.net}
\endgraf{\;}
\endgraf{\;}
\endgraf{Department \, of \, Mathematics Education,\, Kyungpook National University,\,
Taegu \endgraf 702-701, South Korea}
\endgraf {\it Email Address}: {\rm shrim$\@$knu.ac.kr}
\endgraf{\;}
\endgraf{\;}
 }\endaddress
\endthanks
%thanks will be a 1st page footnote
\abstract{The purpose of this paper is to construct new $q$-Euler numbers and polynomials. Finally
we will consider the Witt's type formula associated with these $q$-Euler numbers and polynomials,
and construct $q$-partial zeta functions and $p$-adic $q$-$l$-functions which interpolate
new $q$-Euler numbers and polynomials at negative integers.
}\endabstract
\endtopmatter
\leftheadtext{Taekyun Kim, Min-Soo Kim, Leechae Jang and Seog-Hoon
Rim } \rightheadtext{New $q$-Euler numbers and polynomials}

\document

\head 1. Introduction  \endhead

The constants $E_k$ in the Taylor series expansion
$$\frac2{e^t+1}=\sum_{n=0}^\infty E_n\frac{t^n}{n!}$$
are known as the Euler numbers. The first few are
$1,-\frac12,0,\frac14,-\frac12,\ldots$ and $E_{2k}=0$ for $k=1,2,\ldots.$ Those numbers play an important role
in number theory. For example, the Euler zeta-function essentially equals a Euler numbers at negative integer:
$$\zeta_E(-k)=E_k\quad \text{for } k\geq0,$$
where
$$\zeta_E(s)=\sum_{n=1}^\infty\frac{(-1)^n}{n^s}, \quad s\in \Bbb C \;\,(\text{see \cite{1, 2, 3, 4, 5, 6, 7, 8}}).$$
Recently the $q$-extensions of those Euler numbers and polynomials have been studied by many authors, cf. \cite{1,2,3}.

In \cite{8, 9}, Ozden and Simsek have studied $(h,q)$-extensions of
Euler numbers and polynomials by using $p$-adic $q$-integral on
$\Bbb Z_p.$ From their $(h,q)$-extensions of Euler numbers and
polynomials, they have derived $(h,q)$-extensions of Euler zeta
function and they also gave some interesting relations between their
$(h,q)$-Euler numbers and $(h,q)$-Euler zeta functions, see \cite{8,
9}. Thought this paper $\Bbb Z, \Bbb Z_p, \Bbb Q_p$ and $\Bbb C_p$
will denote the ring of integers, the ring of $p$-adic rational
integers, the the field of $p$-adic rational numbers and completion
of the algebraic closure of $\Bbb Q_p,$ respectively. Let $v_p$ be
the normalized exponential valuation of $\Bbb C_p$ with
$|p|_p=p^{-v_p(p)}=\frac1p.$ When one talks of $q$-extension, $q$ is
variously considered as an indeterminate, a complex number $q\in\Bbb
C,$ or a $p$-adic number $q\in\Bbb C_p.$ If $q\in\Bbb C_p,$ then we
normally assume $|1-q|_p<1.$ If $q\in\Bbb C,$ then we assume that
$|q|<1.$ In this paper we use the following notations:
$$[x]_q=\frac{1-q^x}{1-q}\quad\text{and}\quad [x]_{-q}=\frac{1-(-q)^x}{1+q},\quad \text{cf. \cite{3}}.$$
Let $d$ be a fixed integer, and let
$$\aligned
&X=X_d=\varprojlim_N (\Bbb Z/dp^N\Bbb Z),\quad
X^*=\bigcup\Sb 0<a<dp\\ (a,p)=1\endSb a+dp\Bbb Z_p,\\
&a+dp^N\Bbb Z_p=\{x\in X\mid x\equiv a\pmod{dp^N}\},
\endaligned$$
where $a\in \Bbb Z$ lies in $0\leq a<dp^N$. Let $UD(\Bbb Z_p)$ be the space of uniformly differentiable function on
$\Bbb Z_p.$ For $f\in UD(\Bbb Z_p),$ the $p$-adic $q$-integral was defined by
$$\aligned
I_q(f)&=\int_{\Bbb Z_p}f(x)d\mu_q(x)=\int_Xf(x)d\mu_q(x) \\
&=\lim_{N\rightarrow\infty}\frac1{[dp^N]_q}\sum_{x=0}^{dp^N-1}f(x)q^x\quad\text{for } |1-q|_p<1.
\endaligned$$

In \cite{2, 6} the bosonic integral was considered from a more
physical point of view to the bosonic limit $q\rightarrow1$ as
follows:
$$I_1(f)=\lim_{q\rightarrow1} I_q(f)=\int_{\Bbb Z_p}f(x)d\mu_1(x)=\lim_{N\rightarrow\infty}\frac1{p^N}
\sum_{x=0}^{p^N-1}f(x).$$
Furthermore, we can consider the fermionic integral in contrast to the conventional ``bosonic.''
That is,
$$I_{-1}(f)=\int_{\Bbb Z_p}f(x)d\mu_{-1}(x),\quad\text{see \cite{6}}.\tag1$$
From this, we derive
$$I_{-1}(f_1)+I_{-1}(f)=2f(0),\tag2$$
where $f_1(x)=f(x+1).$ Also we have
$$I_{-1}(f_n)+(-1)^{n-1}I_{-1}(f)=2\sum_{x=0}^{n-1}(-1)^{n-1-x}f(x),$$
where $f_n(x)=f(x+n)$ and $n\in\Bbb Z^+,$ cf. \cite{1,4}.
For $|1-q|_p<1,$ we consider fermionic $p$-adic $q$-integral on $\Bbb Z_p$ which is the $q$-extension of $I_{-1}(f)$
as follows:
$$I_{-q}(f)=\int_{\Bbb Z_p}f(x)d\mu_{-q}(x)
=\lim_{N\rightarrow\infty}\frac1{[dp^N]_{-q}}\sum_{x=0}^{p^N-1}f(x)(-q)^x,
\text{ see \cite{1-11}. }  \tag3$$

By using (1), Ozen and Simsek studied twisted $(h,q)$-Euler numbers
and polynomials and twisted generalized $(h,q)$-Euler numbers
attached to $\chi,$ see \cite{8, 9}.

In this paper, we consider $q$-Euler numbers and polynomials which
are different the $q$-Euler numbers and polynomials of Ozen-Simsek.
Finally, we will give some relations between these $q$-Euler numbers
and polynomials, and construct $q$-partial zeta functions and $p$-adic $q$-$l$-functions which interpolate
new $q$-Euler numbers and polynomials at negative integers.

\head 1. $q$-extensions of Euler numbers and polynomials \endhead

From \cite{6}, we can derive the following formula:
$$qI_{-q}(f_1)+I_{-q}(f)=[2]_qf(0),\tag4$$
where $f_1(x)$ is translation with $f_1(x)=f(x+1).$

If we take $f(x)=e^{tx},$ then we have $f_1(x)=e^{t(x+1)}=e^{tx}e^t.$ From (4), we derive
$$(qe^t+1)I_{-q}(e^{tx})=[2]_q.$$
Hence we obtain
$$I_{-q}(e^{tx})=\int_{\Bbb Z_p}e^{tx}d\mu_{-q}(x)=\frac{[2]_q}{qe^t+1}.\tag5$$
We now define
$$\frac{[2]_q}{qe^t+1}=\sum_{n=0}^\infty E_{n,q}\frac{t^n}{n!}.\tag6$$
By (5) and (6), we see that
$$\int_{\Bbb Z_p}x^nd\mu_{-q}(x)=E_{n,q}.$$
From (4), we also note that
$$\int_{\Bbb Z_p}e^{(x+y)t}d\mu_{-q}(y)=\frac{[2]_q}{qe^t+1}e^{xt}.\tag7$$
In view of (7), we can consider $q$-extension of Euler polynomials as follows:
$$\frac{[2]_q}{qe^t+1}e^{xt}=\sum_{n=0}^\infty E_{n,q}(x)\frac{t^n}{n!}.\tag8$$
By (5), (6), (7) and (8), we obtain the following theorem:

\proclaim{Theorem 1}{\rm(Witt's formula)} For $q\in\Bbb C_p$ with $|1-q|_p<1,$
$$\int_{\Bbb Z_p}x^nd\mu_{-q}(x)=E_{n,q}\quad\text{and}\quad \int_{\Bbb Z_p}(x+y)^nd\mu_{-q}(y)=E_{n,q}(x).$$
\endproclaim

Note that $\lim_{q\rightarrow1}E_{n,q}=E_n$ and $\lim_{q\rightarrow1}E_{n,q}(x)=E_n(x),$ where $E_n$ and $E_{n}(x)$ are
Euler numbers and polynomials.

By Theorem 1, we easily see that
$E_{n,q}(x)=\sum_{k=0}^n\binom nk x^{n-k}E_{k,q}.$
For $n\in\Bbb Z^+,$ let $f_n(x)=f(x+n).$ Then we have
$$q^nI_{-q}(f_n)+(-1)^{n-1}I_{-q}(f)=[2]_q\sum_{l=0}^{n-1}(-1)^{n-l-1}q^lf(l),
\text{ see \cite{2}.}\tag9$$
If $n$ is odd positive integer, we
have
$$q^nI_{-q}(f_n)+I_{-q}(f)=[2]_q\sum_{l=0}^{n-1}(-1)^{l}q^lf(l).\tag10$$
Let $\chi$ be a Dirichlet's character with conductor $d$(=odd)$\in\Bbb Z^+.$ If take $f(x)=\chi(x)e^{tx},$
then we have $f_d(x)=f(x+d)=\chi(x)e^{td}e^{tx}.$
From (10), we derive
$$\int_{X}\chi(x)e^{tx}d\mu_{-q}(x)=\frac{[2]_q\sum_{a=1}^{d}(-1)^aq^a\chi(a)e^{ta}}{q^de^{td}+1}.\tag11$$
In view of (11), we can also consider the generalized Euler numbers attached to $\chi$ as follows:
$$\frac{[2]_q\sum_{a=1}^{d}(-1)^aq^a\chi(a)e^{ta}}{q^de^{td}+1}=\sum_{n=0}^\infty E_{n,\chi,q}\frac{t^n}{n!}.
\tag12$$
From (11) and (12), we derive the following Witt's formula:

\proclaim{Theorem 2} Let $\chi$ be a Dirichlet's character with conductor $d$(=odd)$\in\Bbb Z^+.$ Then we have
$$\int_{X}\chi(x)x^nd\mu_{-q}(x)=E_{n,\chi,q}\tag13$$
for $n\geq0.$
\endproclaim

\head 2. $q$-extension of Euler zeta functions \endhead

For $q\in \Bbb C$ with $|q|<1,$ let
$$F_q(t,x)=\frac{[2]_q}{qe^t+1}e^{tx}=\sum_{n=0}^\infty\frac{E_{n,q}(x)}{n!}t^n\tag14$$
for $|t+\log q|<\pi.$
Then we see that $F_q(t,x)$ is an analytic function on $\Bbb C.$ From (14), we can derive the following expansion:
$$\frac{[2]_q}{qe^t+1}e^{tx}=[2]_q\sum_{n=0}^\infty q^n(-1)^n e^{(n+x)t}.\tag15$$
Thus we have
$$E_{k,q}(x)=\frac{d^k}{dt^k}F_q(t,x)\biggl|_{t=0}=[2]_q\sum_{n=0}^\infty q^n(-1)^n (n+x)^k,\quad k\geq0\tag16$$
and
$$E_{k,q}=\frac{d^k}{dt^k}F_q(t,0)\biggl|_{t=0}=[2]_q\sum_{n=0}^\infty q^n(-1)^n n^k,\quad k\geq0.$$

\proclaim{Definition 3}
For $s\in \Bbb C,$ define
$$\zeta_{q,E}(s,x)=[2]_q\sum_{n=0}^\infty\frac{(-1)^nq^n}{(n+x)^s},
\quad\zeta_{q,E}(s)=[2]_q\sum_{n=1}^\infty\frac{(-1)^nq^n}{n^s}.$$
\endproclaim

Note that $\zeta_{q,E}(s,x)$ and $\zeta_{q,E}(s)$ are analytic functions in the whole complex $s$-plane.

By (15) and (16), we obtain the following:

\proclaim{Theorem 4} Let $n\in\Bbb Z^+\cup\{0\}.$ Then we have
$$\zeta_{q,E}(-n,x)=E_{n,q}(x), \quad \zeta_{q,E}(-n)=E_{n,q}.$$
\endproclaim

Let $\chi$ be a primitive Dirichlet's character with conductor $d$(=odd)$\in\Bbb Z^+.$ Then
the generalized $q$-Euler numbers attached to $\chi$ are defined as
$$F_{q,\chi}(t)=\frac{[2]_q\sum_{a=1}^{d}(-1)^aq^a\chi(a)e^{ta}}{q^de^{td}+1}
=\sum_{n=0}^\infty E_{n,\chi,q} \frac{t^n}{n!},\tag17$$
where $|t+\log q|<\frac \pi d.$

From (17), we note that
$$\aligned
F_{q,\chi}(t)&={[2]_q\sum_{a=1}^{d}(-1)^aq^a\chi(a)e^{ta}}\sum_{l=0}^\infty q^{ld}e^{ldt}(-1)^l \\
&={[2]_q\sum_{n=0}^{\infty}(-1)^nq^n\chi(n)e^{nt}}.
\endaligned\tag18$$

By (17) and (18), we easily see that
$$E_{k,\chi,q}=\frac{d^k}{dt^k}F_{q,\chi}(t)\biggl|_{t=0}=[2]_q\sum_{n=0}^\infty (-1)^nq^n\chi(n)n^k.\tag19$$

\proclaim{Definition 5}
Let $\chi$ be a primitive Dirichlet's character with conductor $d$(=odd)$\in\Bbb Z^+.$ Then
we define the $l_q$-function as follows:
$$l_q(s,\chi)=[2]_q\sum_{n=1}^\infty\frac{(-1)^nq^n\chi(n)}{n^s},\quad s\in\Bbb C.$$
\endproclaim

Note that $l_q(s,\chi)$ is an analytic function in the whole complex $s$-plane.

From (19) and Definition 5, we derive the following:

\proclaim{Theorem 6} For $n\in\Bbb Z^+\cup\{0\},$ we have
$$l_q(-n,\chi)=E_{n,\chi,q}.$$
\endproclaim

Let us consider a partial $q$-zeta function as follows:
$$\aligned H_q(s,a|F)&=[2]_q\sum\Sb{m\equiv a\pmod F}\\m>0\endSb \frac{(-1)^mq^m}{m^s} \\
&=\frac{[2]_q(-1)^aq^a}{F^s}\sum_{n=0}^\infty\frac{(-1)^{nF}q^{nF}}{\left(n+\frac a F\right)^s} \\
&=\frac{(-1)^aq^a}{F^s}\frac{[2]_q}{[2]_{q^F}}\zeta_{q^F,E}(s,\frac a F),
\endaligned\tag20$$
where $F$(=odd) is positive integers with $0<a<F.$ Let $\chi(\neq1)$ be the Dirichlet's character with conductor
$F$(=odd). Then we have
$$l_q(s,\chi)=\sum_{a=1}^F\chi(a)H_q(s,a|F)\tag21$$
for $s\in\Bbb C.$ The function $H_q(s,a|F)$ is an analytic function in whole complex plane. For $n\in\Bbb Z^+,$
we have
$$H_q(-n,a|F)={(-1)^aq^a}{F^n}\frac{[2]_q}{[2]_{q^F}}E_{n,q^F}(\frac aF).\tag22$$
Note that
$$E_{n,q^F}(\frac a F)=\sum_{k=0}^n\binom nk \left(\frac aF\right)^{n-k}E_{k,q^F}.\tag23$$
By using (22) and (23) we have
$$H_q(-n,a|F)={(-1)^aq^a}{F^n}\frac{[2]_q}{[2]_{q^F}}\sum_{k=0}^n\binom nk \left(\frac aF\right)^{n-k}E_{k,q^F}.\tag24$$
We now modify a partial $q$-zeta function as follows:
$$H_q(s,a|F)={(-1)^aq^a}{a^{-s}}
\frac{[2]_q}{[2]_{q^F}}\sum_{k=0}^\infty\binom {-s}k \left(\frac Fa\right)^{k}E_{k,q^F}\tag24\text{$^\prime$}$$
for $s\in\Bbb C.$ From (21) and (24$^\prime$), we note that
$$l_q(s,\chi)=\frac{[2]_q}{[2]_{q^F}}\sum_{a=1}^F(-1)^a\chi(a){q^a}{a^{-s}}\sum_{k=0}^\infty\binom {-s}k \left(\frac Fa\right)^{k}
E_{k,q^F}\tag25$$
for $s\in\Bbb C.$ By (25), we easily see that
$$l_q(s,\chi)=\frac{[2]_q}{[2]_{q^F}}\sum_{a=1}^F(-1)^a\chi(a){q^a}\left\{a^{-s}+a^{-s}\sum_{k=1}^\infty\binom {-s}k
E_{k,q^F}\left(\frac Fa\right)^{k}\right\}.\tag26$$
From the Taylor expansion at $s=0,$ we have
$$l_q(0,\chi)=\frac{[2]_q}{[2]_{q^F}}\sum_{a=1}^F(-1)^aq^a\chi(a),$$
and
$$l_q(1,\chi)=\frac{[2]_q}{[2]_{q^F}}\sum_{a=1}^F\frac{(-1)^aq^a}{a}\chi(a)\left\{1+
\sum_{k=1}^\infty(-1)^k E_{k,q^F}\left(\frac Fa\right)^{k}\right\}.$$

\head 3. $p$-adic interpolating function for $q$-Euler numbers\endhead

We shall consider the $p$-adic analogue of the $l_q$-functions which are introduced in the previous section.
Throughout this section we assume that $p$ is an odd prime.
Let $\omega$ be denoted as the Teichm\"uller character having conductor $p.$ For an arbitrary character $\chi,$
let $\chi_n=\chi\omega^{-n},$ where $n\in \Bbb Z,$ in sense of the product of characters.
Let $\langle a\rangle=\omega^{-1}(a)a=\frac{a}{\omega(a)}.$

Let $\chi$ be the Dirichlet's character with conductor $d$(=odd) and let $F$ be a positive integral multiple
of $p$ and $d.$ Now, we define the $p$-adic $l_q$-functions as follows:
$$l_{p,q}(s,\chi)=\frac{[2]_q}{[2]_{q^F}}
\sum\Sb a=1\\ p\nmid a\endSb^F(-1)^a\chi(a){q^a}\langle a\rangle^{-s}\sum_{k=0}^\infty\binom {-s}k \left(\frac Fa\right)^{k}
E_{k,q^F}.\tag27$$
Then $l_q(s,\chi)$ is an analytic function in $D=\{s\in\Bbb C_p\mid |s|_p<p^{1-\frac1{p-1}}\}$ since
$\langle a\rangle^{-s}$ and $\sum_{k=0}^\infty\binom {-s}k \left(\frac Fa\right)^{k}
E_{k,q^F}$ are analytic functions in $D,$ cf. \cite{5,6,7,9}.

We set
$$H_{p,q}(s,a|F)={(-1)^aq^a}\langle a\rangle^{-s}
\frac{[2]_q}{[2]_{q^F}}\sum_{k=0}^\infty\binom {-s}k \left(\frac Fa\right)^{k}E_{k,q^F}.\tag28$$
Thus, by (24), note that
$$H_{p,q}(-n,a|F)=\omega^{-n}(a)H_q(-n,a|F).\tag29$$
for $n\in\Bbb Z^+.$ We also consider the $p$-adic analytic
function which interpolates $q$-Euler number at negative integer
as follows:
$$l_{p,q}(s,\chi)=\sum\Sb a=1\\ p\nmid a\endSb^F\chi(a)H_{p,q}(s,a|F).\tag30$$
For $d$(=odd)$\in\Bbb Z^+,$ by (14) and (17), we note that
$$\aligned
\frac1{[2]_q}F_{q,\chi}(t)&=\frac1{[2]_{q^d}}\sum_{a=1}^{d}(-1)^aq^a\chi(a)F_{q^d}(dt,\frac ad) \\
&=\sum_{n=0}^\infty\left[d^n\frac1{[2]_{q^d}}\sum_{a=1}^{d}(-1)^aq^a\chi(a)E_{n,q^d}(\frac ad)\right]\frac{t^n}{n!}.
\endaligned\tag31$$
Then we have
$$\frac1{[2]_{q}}E_{n,\chi,q}=d^n\frac1{[2]_{q^d}}\sum_{a=1}^{d}(-1)^aq^a\chi(a)E_{n,q^d}(\frac ad).\tag32$$
In particular, if $F=dp,$ then we have
$$\frac1{[2]_{q}}E_{n,\chi,q}=F^n\frac1{[2]_{q^F}}\sum_{a=1}^{F}(-1)^aq^a\chi(a)E_{n,q^F}(\frac aF).\tag32\text{$^\prime$}$$
Indeed, by (31), it is sufficient to show that
$$\aligned
\frac1{[2]_{q^F}}\sum_{a=1}^{F}(-1)^aq^a\chi(a)F_{q^F}(Ft,\frac aF)&=\sum_{a=1}^{F}(-1)^aq^a\chi(a)
\frac{e^{ta}}{q^Fe^{Ft}+1} \\
&=\sum_{b=1}^d\sum_{c=0}^{p-1}(-1)^{b+cd}\chi({b+cd})q^{b+cd}\frac{e^{t({b+cd})}}{q^Fe^{Ft}+1}\\
&=\sum_{b=1}^d(-1)^bq^b\chi(b)\frac{e^{bt}}{q^Fe^{Ft}+1}\sum_{c=0}^{p-1}(-1)^c(q^de^{dt})^c \\
&=\frac1{[2]_{q^d}}\sum_{b=1}^d(-1)^bq^b\chi(b)F_{q^d}(dt,\frac bd).
\endaligned$$
Thus, since $\chi_n=\chi\omega^{-n},$ by (29) and (32\text{$^\prime$}), we obtain
$$\aligned
l_{p,q}(-n,\chi)&=\sum\Sb a=1\\ p\nmid a\endSb^F\chi(a)H_{p,q}(-n,a|F) \\
&=\sum\Sb a=1\\ p\nmid a\endSb^F\chi_n(a)H_{q}(-n,a|F) \\
&=F^n\frac{[2]_q}{[2]_{q^F}}\sum\Sb a=1\\ p\nmid a\endSb^F\chi_n(a)(-1)^aq^aE_{n,q^F}(\frac aF)\\
&=F^n\frac{[2]_q}{[2]_{q^F}}\sum_{a=1}^F\chi_n(a)(-1)^aq^aE_{n,q^F}(\frac aF) \\
&\qquad-F^n\frac{[2]_q}{[2]_{q^F}}\sum_{a=1}^{\frac Fp}\chi_n(pa)(-1)^{pa}q^{pa}E_{n,q^{F}}(\frac {pa}F).
\endaligned\tag33$$
From Theorem 2 and (32\text{$^\prime$}), we see that
$$\frac1{[2]_q}E_{n,\chi,q}=\frac1{[2]_{q^F}}F^n\sum_{a=1}^{F}(-1)^aq^a\chi(a)E_{n,q^F}(\frac aF)
=\frac1{[2]_q}\int_{X}\chi(x)x^nd\mu_{-q}(x).\tag34$$
and
$$\frac1{[2]_q}E_{n,\chi,q^p}=\left(\frac Fp\right)^n\frac1{[2]_{q^F}}
\sum_{a=1}^{\frac Fp}(-1)^a(q^{p})^a\chi(a)E_{n,(q^{p})^{\frac Fp}}
(\frac a{\frac Fp}).\tag35$$
Also, by (34) and (35), we have
$$\aligned\int_{pX}\chi(x)x^nd\mu_{-q}(x)&=\int_{X}\chi(px)(px)^nd\mu_{-q}(px) \\
&=\chi(p)p^n\frac{[2]_q}{[2]_{q^p}}\int_{X}\chi(x)x^nd\mu_{-q^p}(x) \\
&=\chi(p)p^n\frac{[2]_q}{[2]_{q^p}}F^n\sum_{a=1}^{F}(-1)^a(q^{p})^a\chi(a)E_{n,(q^{p})^F}(\frac aF)\\
&=\chi(p)p^n\frac{[2]_q}{[2]_{q^p}}E_{n,\chi,q^p},
\endaligned$$
since $d\mu_{-q}(px)=\frac{[2]_q}{[2]_{q^p}}d\mu_{-q^p}(x).$ Therefore, we obtain the following theorem:

\proclaim{Theorem 7}
Let $F$(=odd) be a positive integral multiple of $p$ and $d=(d_\chi),$ and let
$$l_{p,q}(s,\chi)=\frac{[2]_q}{[2]_{q^F}}\sum\Sb a=1\\ p\nmid a\endSb^F(-1)^a\chi(a){q^a}\langle a\rangle^{-s}\sum_{k=0}^\infty\binom {-s}k
\left(\frac Fa\right)^{k}E_{k,q^F}.$$
Then we have
\roster
\item"(a)" $l_{p,q}(s,\chi)$ analytic in $D=\{s\in\Bbb C_p\mid |s|_p<p^{1-\frac1{p-1}}\}.$
\item"(b)" $l_{p,q}(-n,\chi)=E_{n,\chi_n,q}-p^n\chi_n(p)\frac{[2]_q}{[2]_{q^p}}E_{n,\chi_n,q^p}$ for $n\in\Bbb Z^+.$
\item"(c)" For $n\in\Bbb Z^+,$
$$l_{p,q}(-n,\chi)=\int_{X^*}\chi_n(x)x^nd\mu_{-q}(x).$$
\endroster
\endproclaim

\proclaim{Corollary 8} Let $F$(=odd) be a positive integral multiple
of $p$ and $d=(d_\chi),$ and let
$$l_{p}(s,\chi)=\sum\Sb a=1\\ p\nmid a\endSb^F(-1)^a\chi(a)\langle a\rangle^{-s}\sum_{k=0}^\infty\binom {-s}k
\left(\frac Fa\right)^{k}E_{k}, \text{ see [10]}.$$ Then we have
\roster
\item"(a)" $l_{p}(s,\chi)$ analytic in $D=\{s\in\Bbb C_p\mid |s|_p<p^{1-\frac1{p-1}}\}.$
\item"(b)" $l_{p,q}(-n,\chi)=E_{n,\chi_n}-p^n\chi_n(p)E_{n,\chi_n}$ for $n\in\Bbb Z^+.$
\item"(c)" $l_{p}(s,\chi)=\int_{X^*}\chi_n(x)x^{-s}d\mu_{-1}(x).$ Observe that for $n\in\Bbb Z^+,$
$$l_{p}(-n,\chi)=\int_{X^*}\chi_n(x)x^nd\mu_{-1}(x).$$
\endroster
\endproclaim

\remark{Remark 9}
In the recent paper (see \cite{9}), Ozden and Simsek have studied the $(h,q)$-extension of twisted Euler numbers.
However, these $(h,q)$-extension of twisted Euler numbers and generating function do not seen to be natural ones;
in particular, these numbers cannot be represented as a nice Witt's type formula for the $p$-adic
invariant integral on $\Bbb Z_p$ and the generating function does not seems to be simple and useful for deriving
many interesting identities related to the extension of Euler numbers. By this reason, we consider new $q$-extensions
of Euler numbers and polynomials which are different. Our $q$-extensions
of Euler numbers and polynomials to treat in this paper can be represented by $p$-adic $q$-fermionic integral on
$\Bbb Z_p$ and this integral reprentation also can consider as Witt's type formula for $q$-extensions
of Euler numbers and polynomials.
\endremark

\enddocument